\documentclass[11pt,a4paper]{article}
\pdfoutput=1 
\usepackage[margin=1in]{geometry}
\usepackage{amsmath,amsthm,amssymb,enumerate, bbm ,graphicx,color,caption,upgreek, float, tikz, wasysym, subcaption,booktabs,longtable, appendix,graphics, pdfpages,rotating,mathtools,textcomp, array, blkarray}

\DeclarePairedDelimiter\ceil{\lceil}{\rceil}
\DeclarePairedDelimiter\floor{\lfloor}{\rfloor}
\usepackage[pdftex]{hyperref}

\definecolor{blauw}{RGB}{61,158,255}
\definecolor{donkerblauw}{RGB}{0,0,255}
\definecolor{donkergroen}{RGB}{46,148,0}
\definecolor{donkerrood}{RGB}{204,0,0}

\makeatletter 
\newcommand\mynobreakpar{\par\nobreak\@afterheading} 
\makeatother

\makeatletter
\let\@fnsymbol\@arabic
\makeatother

\newcommand{\N}{\mathbb{N}}

\usepackage[english]{babel}

\newtheorem{theorem}{Theorem}[section]

\newtheorem{claim}[theorem]{Claim}
\newtheorem{proposition}[theorem]{Proposition}
\newtheorem{corollary}[theorem]{Corollary}

\theoremstyle{definition}
\newtheorem{defn}{Definition}[section]
\newtheorem{assumption}{Assumption}[section]
\newtheorem{examp}{Example}[section]
\newtheorem*{examp*}{Example}

\newtheorem{remark}{Remark}[section]

\theoremstyle{plain}

\newfloat{Algorithm}{!hbt}{alg}

\newcounter{thm}[section]

\title{New nonbinary code bounds based on divisibility arguments} \date{}
\author{Sven Polak\thanks{Korteweg-De Vries Institute for Mathematics, University of Amsterdam. E-mail: \href{mailto:s.c.polak@uva.nl}{s.c.polak@uva.nl}. The research leading to these
results has received funding from the European Research Council under the European Union’s Seventh Framework Programme (FP7/2007-2013) / ERC grant agreement \textnumero 339109.}}
\begin{document}
\maketitle
\setcounter{footnote}{1}

\noindent \textbf{Abstract.} For~$q,n,d \in \N$, let~$A_q(n,d)$ be the maximum size of a code~$C \subseteq [q]^n$ with minimum distance at least~$d$. We give a divisibility argument resulting in the new upper bounds~$A_5(8,6) \leq 65$, $A_4(11,8)\leq 60$ and~$A_3(16,11) \leq 29$. These in turn imply the new upper bounds~$A_5(9,6) \leq 325$,~$A_5(10,6) \leq 1625$,~$A_5(11,6) \leq 8125$ and~$A_4(12,8) \leq 240$.

Furthermore, we prove that for~$\mu,q \in \N$,  there is a 1-1-correspondence between symmetric~$(\mu,q)$-nets (which are certain designs) and codes~$C \subseteq [q]^{\mu q}$ of size~$\mu q^2$ with minimum distance at least~$\mu q - \mu$. We derive the new upper bounds~$A_4(9,6) \leq 120$ and~$A_4(10,6) \leq 480$ from these `symmetric net' codes. 

\,$\phantom{0}$


\section{Introduction}
\noindent For any~$m \in \mathbb{N}$, we write~$[m]:=\{1,\ldots,m\}$. Fix~$n ,q \in \N$. A \emph{word} is an element~$v \in [q]^n$. So~$[q]$ serves as the alphabet. (If you prefer~$\{0,1,\ldots,q-1\}$ as alphabet, take the letters mod~$q$.)  For two words~$u,v \in [q]^n$, their \emph{(Hamming) distance}~$d_H(u,v)$ is the number of indices~$i$ with~$u_i \neq v_i$. A \emph{code} is a subset of~$ [q]^n$. For any code~$C \subseteq  [q]^n$, the minimum distance~$d_{\text{min}}(C)$ of~$C$ is the minimum distance between any two distinct code words in~$C$. For~$d \in \N$, an~$(n,d)_q$\emph{-code} is a set~$C\subseteq [q]^n$ that satisfies~$d_{\text{min}}(C)\geq d$. Define
\begin{align}\label{aqnd}
    A_q(n,d) :=  \max \{ | C| \,\, | \,\, C  \text{ is an $(n,d)_q$-code} \}.
\end{align}
Computing~$A_q(n,d)$ and finding upper and lower bounds for it is a long-standing research interest in combinatorial coding theory (cf.\ MacWilliams and Sloane~\cite{sloane}). In this paper we find new upper bounds on~$A_q(n,d)$ (for some~$q,n,d$), based on a \emph{divisibility}-argument. In some cases, it will sharpen a combination of the following two well-known upper bounds on~$A_q(n,d)$. Fix~$q,n,d \in \N$. Then
\begin{align} \label{elementarybounds}
qd >(q-1)n\,\,\, \Longrightarrow\,\,\, A_q(n,d) \leq \frac{qd}{qd-n(q-1)}.
\end{align}
This is the~$q$-ary \emph{Plotkin} bound. Moreover,
\begin{align} \label{elementarybounds2}
A_q(n,d) \leq q \cdot A_q(n-1,d). 
\end{align}

\noindent A proof of these statements can be found in~$\cite{sloane}$.  Plotkin's bound can be proved by comparing the leftmost and rightmost terms in~$(\ref{detruc})$ below. The second bound follows from the observation that in a~$(n,d)_q$-code any symbol can occur at most~$A_q(n-1,d)$ times at the first position.

We view an~$(n,d)_q$-code~$C$ of size~$M$ as an~$M \times n$ matrix with the words as rows. Two codes~$C,D \subseteq [q]^n$ are \emph{equivalent} (or \emph{isomorphic}) if~$D$ can be obtained from~$C$ by first permuting the~$n$ columns of~$C$ and subsequently applying to each column a permutation of the~$q$ symbols in~$[q]$ (we will write `\emph{renumbering a column}' instead of `applying a permutation to the symbols in a column').

\begin{table}[H]
\begin{center}
    \begin{tabular}{| l || r |r ||>{\bfseries}r|}
    \hline
    $A_q(n,d)$ &  lower bound & upper bound & new upper bound \\
         &$\cite{4ary,5ary,plotkin}$   & $\cite{4ary,5ary}$  &\\  \hline 
    $A_5(8,6)$ & 50 & 75 &  65 \\ 
    $A_5(9,6)$ & 135 & 375  & 325 
    \\ $A_5(10,6)$ & 625 & 1855  & 1625 \\    
        $A_5(11,6)$ & 3125 & 8840  & 8125 \\   
    \hline
     $A_4(9,6)$ &  64  &  128     & 120 \\  
          $A_4(10,6)$ & 256  & 496   &480  \\      
    $A_4(11,8)$ & 48 & 64  & 60  \\    
    $A_4(12,8)$ & 128 & 242   & 240  \\    
     \hline 
    $A_3(16,11)$ & 18 & 30    & 29  \\   \hline  
    \end{tabular}
\end{center}
  \caption{An overview of the results obtained and discussed in this paper. All previous lower and upper bounds are taken from references~$\cite{4ary,5ary}$, except for the lower bounds~$A_5(8,6)\geq 50$ and~$A_4(11,8) \geq 48$.\protect\footnotemark\, These follow from the exact values~$A_5(10,8)=50$ and~$A_4(12,9)=48$ ($\cite{plotkin}$). For updated tables with all most recent code bounds, we refer to~$\cite{brouwertable}$. }
\end{table} 
\footnotetext{In~$\cite{4ary,5ary}$, the lower bounds~$A_5(8,6)\geq 45$ and~$A_4(11,8) \geq 34$ are given.}

 If an~$(n,d)_q$-code~$C$ is given, then for~$j=1,\ldots,n$, let~$c_{\alpha,j}$ denote the number of times symbol~$\alpha \in [q]$ appears in column~$j$ of~$C$. For any two words~$u,v \in [q]^n$, we define~$g(u,v):=n-d_H(u,v)$.  In our divisibility arguments, we will use the following observations (which are well known and often used in coding theory and combinatorics). 
 
\begin{proposition} 
If~$C$ is an~$(n,d)_q$-code of size~$M$, then
\begin{align}\label{detruc}
\binom{M}{2}(n-d) \geq \sum_{\substack{\{u,v\} \subseteq C\\ u \neq v}} g(u,v) = \sum_{j=1}^n \sum_{\alpha \in [q]} \binom{c_{\alpha,j}}{2} \geq n \cdot \left( (q-r)\binom{m}{2} +r \binom{m-1}{2}  \right),
\end{align}
where~$m := \ceil{M/q}$ and~$ r:=qm-M$, so that~$M=qm-r$ and~$0 \leq r < q$. Moreover, writing~$L$ and~$R$ for the leftmost term and the rightmost term in~$(\ref{detruc})$, respectively, we have 
\begin{align} \label{boundnd} 
| \{ \{u,v\} \subseteq C \,\, | \,\, u \neq v, \,\, d_H(u,v)\neq d \}| \leq L - R,  
\end{align}
i.e., the number of pairs of distinct words~$\{u,v\} \subseteq C$ with distance \emph{unequal} to~$d$ is at most the leftmost term minus the rightmost term in~$(\ref{detruc})$.
\end{proposition}
\proof 
The first inequality in~$(\ref{detruc})$ holds because~$n-d \geq g(u,v)$ for all~$u,v \in C$. The equality is obtained by counting the number of equal pairs of entries in the same columns of~$C$ in two ways. The second inequality follows from the (strict) convexity of the binomial coefficient~$F(x):=x(x-1)/2$. Fixing a column~$j$, the quantity~$\sum_{\alpha \in [q]} F(c_{\alpha,j})$, under the condition that~$\sum_{\alpha \in [q]} c_{\alpha,j} = M$, is minimal if the~$c_{\alpha,j}$ are as equally divided as possible, i.e., if~$c_{\alpha,j} \in \{ \ceil{M/q}, \floor{M/q}\}$ for all~$\alpha \in [q]$. The desired inequality follows.

To prove the second assertion, note that it follows from~$(\ref{detruc})$ that~$\sum_{\{u,v\} \subseteq C, \, u \neq v} g(u,v) \geq R$, so
\begin{align}\label{detruc2}
| \{ \{u,v\} \subseteq C \,\, | \,\, u \neq v, \,\, d_H(u,v)\neq d \}|  &\leq \sum_{\substack{\{u,v\} \subseteq C\\ u \neq v}} (n-d-g(u,v)) 
\\&\leq  \binom{M}{2}(n-d) - R = L-R.\notag \qedhere
\end{align}
\endproof

\begin{corollary}\label{cor} If, for some~$q, n, d$ and~$M$, the left hand side equals the right hand side in~$(\ref{detruc})$, then for any $(n,d)_q$-code~$C$ of size~$M$, 
\begin{enumerate}[(i)] 
\item \label{1} $g(u,v)=n-d$ for all~$u,v \in C$ with~$u \neq v$, i.e.,~$C$ is \emph{equidistant}, and
\item \label{2} for each column~$C_j$ of~$C$, there are $q-r$ symbols in~$[q]$ that occur~$m$ times in~$C_j$ and~$r$ symbols in~$[q]$ that occur~$m-1$ times in~$C_j$. 
\end{enumerate}
\end{corollary}
\noindent In the next sections we will use~(\ref{1}),~(\ref{2}) and the bound in~$(\ref{boundnd})$ to give (for some~$q,n,d$) new upper bounds on~$A_q(n,d)$, based on divisibility arguments. Furthermore, in Section~$\ref{symsec}$, we will prove that, for~$\mu,q \in \N$, there is a 1-1-correspondence between symmetric~$(\mu,q)$-nets (which are certain designs) and $(n,d)_q=(\mu q, \mu q - \mu)_q$-codes~$C$ with~$|C|=\mu q^2 $. We derive some new upper bounds from these `symmetric net' codes.

\section{The divisibility argument}
In this section, we describe the divisibility argument and illustrate it by an example. Next, we show how the divisibility argument can be applied to obtain upper bounds on~$A_q(n,d)$ for certain~$q,n,d$. In subsequent sections, we will see how we can improve upon these bounds for certain fixed~$q,n,d$. We will use the following notation. 

\begin{defn}[$k$-block]
Let~$C$ be an~$(n,d)_q$-code in which a symbol~$\alpha \in [q]$ is contained exactly~$k$ times in column~$j$. The~$k \times n$ matrix~$B$ formed by the $k$ rows of~$C$ that have symbol~$\alpha$ in column~$j$ is called a~($k$-)\emph{block} (for column~$j$). In that case, columns~$[n]\setminus \{j\}$ of~$B$ form an~$(n-1,d)_q$-code of size~$k$. 
\end{defn}

\noindent At the heart of the divisibility arguments that will be used throughout this paper lies the following observation.
\begin{proposition}[Divisibility argument] \label{parity}
Suppose that~$C$ is an~$(n,d)_q$-code and that~$B$ is a block in~$C$ (for some column~$j$) containing every symbol exactly~$m$ times in every column except for column~$j$. If~$n-d$ does not divide~$m(n-1)$, then for each~$u \in C \setminus B$ there is a word~$v \in B$ with~$d_H(u,v) \notin \{d,n\}$.
\end{proposition}
\proof 
Let~$u \in C \setminus B$. We renumber the symbols in each column such that~$u$ is $\mathbf{1}:=1\ldots 1$, the all-ones word. The total number of~1's in~$B$ is~$m(n-1)$ (as the block~$B$ does not contain~1's in column~$j$ since~$u\notin B$ and~$B$ consists of all words in~$C$ that have the same symbol in column~$j$). Since~$n-d$ does not divide~$m(n-1)$, there must be a word~$v \in B$ that contains a number of~1's not divisible by~$n-d$. In particular, the number of~$1$'s in~$v$ is different from~$0$ and~$n-d$. So~$d_H(u,v) \notin \{d,n\}$.
\endproof 

\begin{examp}\label{586}
We apply Proposition~$\ref{parity}$ to the case~$(n,d)_q=(8,6)_5$. The best known upper bound\footnote{The Delsarte bound~$\cite{delsarte}$ on~$A_5(8,6)$, the bound based on Theorem~$\ref{elementarybounds}$, and the semidefinite programming bound based on quadruples of code words~\cite{onsartikel} all are equal to 75.} is~$A_5(8,6) \leq 75$, which can be derived from~$(\ref{elementarybounds})$ and~$(\ref{elementarybounds2})$, as the Plotkin bound yields~$A_5(7,6) \leq 15$ and hence~$A_5(8,6) \leq 5 \cdot 15 = 75$. Since, for~$(n,d)_q=(7,6)_5$ and~$M=15$, the left hand side equals the right hand side in~$(\ref{detruc})$, any $(7,6)_5$-code~$D$ of size~$15$ is equidistant and each symbol appears exactly~$m=3$ times in every column of~$D$. Note~$2=n-d \nmid m(n-1)=21$.

Suppose there exists a~$(8,6)_5$-code~$C$ of size~$75$. As~$A_5(7,6) \leq 15$, for each column,~$C$ is divided into five 15-blocks. Let~$B$ be a~$15$-block for the~$j$th column and let~$u \in C \setminus B$. Note that the other columns of~$B$ contain each symbol~$3$ times, and~$3(n-1)=3\cdot 7 =21$ is not divisible by~$n-d=2$. So by Proposition~\ref{parity}, there must be a word~$v \in  B$ with~$d_H(u,v) \notin \{6,8\}$.
 
However, since all~$(7,6)_5$-codes of size~15 are equidistant, all distances in~$C$ belong to~$\{6,8\}$: either two words are contained together in some~$15$-block (hence their distance is~$6$) or there is no column for which the two words are contained in a~$15$-block (hence their distance is~$8$). This implies that an~$(8,6)_5$-code~$C$ of size~$75$ cannot  exist. Hence~$A_5(8,6) \leq 74$. Theorem~$\ref{importantth}$ and Corollary~$\ref{1mod4}$ below will imply that~$A_5(8,6)\leq 70$ and in Section~$\ref{sec586}$ we will show that, with some computer assistance, the bound can be pushed down to~$A_5(8,6) \leq 65$.
\end{examp}

\noindent To exploit the idea of Proposition~$\ref{parity}$, we will count the number of so-called irregular pairs of words occuring in a code.

\begin{defn}[Irregular pair]
Let~$C$ be an~$(n,d)_q$-code and~$u,v \in C$ with~$u\neq v$. If~$d_H(u,v) \notin \{d,n\}$, we call~$\{u,v\}$ an \emph{irregular pair}. 
\end{defn}

\noindent For any code~$C \subseteq  [q]^n$, we write
\begin{align}\label{X}
X:= \text{ the set of irregular pairs~$\{u,v\}$ for~$u,v \in C$}.
\end{align}

\noindent Using Proposition~$\ref{parity}$, we can for some cases derive a lower bound on~$|X|$. If we can also compute an upper bound on~$|X|$ that is smaller than the lower bound, we derive that the code~$C$ cannot exist. The proof of the next theorem uses this idea. For fixed~$q,n,d,m \in \N$ with~$q \geq 2$, define the following quadratic polynomial in~$r$:
\begin{align} 
\phi(r) := n(n-1-d)(r-1)r - (q-r+1)(mq(q+r-2)-2r). 
\end{align}
\begin{theorem} \label{importantth}
Suppose that~$q\geq 2$, that~$m:=d/(qd-(n-1)(q-1))$ is a positive integer, and that~$n-d$ does not divide~$m(n-1)$. If~$r \in \{1,\ldots,q-1\}$ with~$\phi(r)<0$, then~$A_q(n,d) < mq^2 -r$. \end{theorem}
\proof 
 By Plotkin's bound~$(\ref{elementarybounds})$ we have
\begin{align} \label{part}
    A_q(n-1,d) \leq mq. 
\end{align}
Let~$D$ be an~$(n-1,d)_q$-code of size~$mq-t$ with~$t<q$. Note that~$d=m(n-1)(q-1)/(mq-1)$. Then the right-hand side in~$(\ref{boundnd})$ (taking~$C:= D)$  is equal to~$(n-1)(m-1)t(t-1)/(2mq-2) = (n-1-d)\binom{t}{2}$. Hence
\begin{align} \label{most}
    \text{$D$ contains at most~$(n-1-d)\binom{t}{2}$ pairs of words with distance~$\neq d$}.
\end{align}
Therefore, all~$(n-1,d)_q$-codes~$D$ of size~$mq$ are equidistant (then~$t=0$) and each symbol occurs~$m$ times in every column of~$D$.

 Now let~$C$ be an~$(n,d)_q$-code of size~$M:=mq^2-r$ with~$r \in \{1,\ldots,q-1\}$. Consider an~$mq$-block~$B$ for some column of~$C$. As~$n-d$ does not divide~$m(n-1)$, by Proposition~$\ref{parity}$ we know 
\begin{align}\label{obs}
\text{if~$u \in C\setminus B$, then there exists~$v \in B$ with~$d_H(u,v) \notin \{d,n\}$.}
\end{align}
Let~$B_1,\ldots, B_{s}$ be $mq$-blocks~in~$C$ for some fixed column. Since~$|C|=mq^2-r$, the number of $mq$-blocks for any fixed column is at least~$q-r$ (so we can take~$s=q-r$). Then, with~$(\ref{obs})$, one obtains a lower bound on the number~$|X|$ of irregular pairs in~$C$. Every pair~$\{B_{i},B_{k}\}$ of $mq$-blocks gives rise to~$mq$ irregular pairs: for each word~$u \in B_{i}$, there is a word~$v \in B_k$ such that~$\{u,v\}\in X$. This implies that in~$\cup_{i=1}^s B_i \subseteq C$ there are at least~$\binom{s}{2}mq$ irregular pairs. Moreover, for each word~$u$ in~$C\setminus \cup_{i=1}^s B_i$ (there are~$M-mq\cdot s$ of such words) there is, for each~$i=1,\ldots,s$, a word~$v_i \in B_i$  with~$\{u,v_i\} \in X$. This gives an additional number of at least~$(M-mqs)s$ irregular pairs in~$C$. Hence:
\begin{align} \label{lowerb}
|X| &\geq \binom{s}{2} mq + (M-mqs)s \notag  
\\&= \mbox{$\frac{1}{2}$}s (mq(2q-s-1)-2r) =: l(s).
\end{align} 

On the other hand, note that the~$i$th block for the~$j$th column has size~$mq-r_{i,j}$ for some integer~$r_{i,j}\geq 0$ by~$(\ref{part})$, where~$\sum_{i=1}^q r_{i,j}=r\leq q-1$ (hence each~$r_{i,j} <q$). So by~$(\ref{most})$, the number of irregular pairs in~$C$ that have the same entry in column~$j$ is at most
\begin{align} 
(n-1-d)\sum_{i=1}^q \binom{r_{i,j}}{2}.
\end{align} 
 As each irregular pair~$\{u,v\}$ has~$u_j=v_j$ for at least one column~$j$, we conclude  
\begin{align} \label{upperb}
|X| \leq  (n-1-d)\sum_{j=1}^n \sum_{i=1}^q \binom{r_{i,j}}{2} \leq n (n-1-d)\binom{r}{2}.
\end{align} 
Here the last inequality follows by convexity of the binomial function, since (for fixed~$j$) the sum~$\sum_{i=1}^q \binom{r_{i,j}}{2}$ under the condition that~$\sum_{i=1}^q r_{i,j}=r$  is maximal if one of the~$r_{i,j}$ is equal to~$r$ and the others are equal to~$0$.

 If each~$r_{i,j}\in \{0,1\}$, then~$|X|=0$ by~$(\ref{upperb})$. As~$q-r\geq 1$, there is at least one~$mq$-block for any fixed column, so~$|X|\geq 1$ by~$(\ref{obs})$, which is not possible. Hence we can assume that~$r_{i,j} \geq 2$ for some~$i,j$ (this also implies~$A_q(n,d) \leq mq^2-2$). Then the number~$s$ of~$mq$-blocks for column~$j$ satisfies~$s \geq q-r+1$. This gives by~$(\ref{lowerb})$ and~$(\ref{upperb})$ that
 \begin{align} \label{newineq2}
l(q-r+1) \leq |X| \leq (n-1-d) \binom{r}{2}.
\end{align}
Subtracting the left hand side from the right hand side in~$(\ref{newineq2})$ yields~$\phi(r)/2 \geq 0$, i.e.,~$\phi(r) \geq 0$. So if~$\phi(r)< 0$, then~$A_q(n,d) < mq^2-r$, as was needed to prove.
\endproof

\noindent We give two interesting applications of Theorem~$\ref{importantth}$.

\begin{corollary} \label{1mod4}
If~$q \equiv 1 \pmod{4}$ and~$q\neq 1$, then
 \begin{align}\label{q+3}
     A_q(q+3,q+1) \leq \mbox{$\frac{1}{2}$}q^2(q+1)-q =\mbox{$\frac{1}{2}$}(q-1)q (q+2). 
 \end{align}
\end{corollary}
\proof 
Apply Theorem~$\ref{importantth}$ to~$n=q+3$,~$d=q+1$ and~$r=q-1$. Then~$m=(q+1)/2 \in \N$ and~$n-d=2$ does not divide~$m(n-1)=(q+1)(q+2)/2$, as~$q\equiv 1 \pmod{4}$. Furthermore, $\phi(q-1)=-(q^3-q^2-2)<0$.  Hence~$A_q(q+3,q+1) < q^2(q+1)/2-(q-1)$.
\endproof 

\noindent Applying Corollary~$\ref{1mod4}$ to~$q=5$ gives~$A_5(8,6) \leq 70$. In Section~$\ref{sec586}$ we will improve this to~$A_5(8,6) \leq 65$.

\begin{remark}
Note that for bound~$(\ref{q+3})$ to hold it is necessary that~$q\equiv 1 \pmod{4}$. If~$q \equiv 3 \pmod{4}$ the statement does \emph{not} hold in general. For example,~$A_3(6,4)= 18$ (see~$\cite{brouwertable}$), which is larger than bound~$(\ref{q+3})$. 
\end{remark}

\noindent Theorem~$\ref{importantth}$  also gives an upper bound on~$A_q(n,d)=A_q(kq+k+q,kq)$, where~$q\geq 2$ and~$k$ does not divide~$q(q+1)$ (which is useful for~$k <q-1$; for~$k \geq q+1$ the Plotkin bound gives a better bound). One new upper bound for such~$q,n,d$ is obtained:

\begin{proposition}\label{60} $A_4(11,8)\leq 60$. 
\end{proposition}
\proof 
This follows from Theorem~$\ref{importantth}$ with~$q=4$, $n=11$, $d=8$ and~$r=3$. Then~$m=4 \in \N$, and~$n-d=3$ does not divide~$m(n-1)=40$. Moreover,~$\phi(3)=-16 <0$. Therefore~$A_4(11,8) < 61$.
\endproof  

\noindent This implies the following bound, which is also new:

\begin{corollary} 
$A_4(12,8) \leq 240$.
\end{corollary}
\proof 
By Proposition~$\ref{60}$ and~$(\ref{elementarybounds2})$.
\endproof

\section{Kirkman triple systems and \texorpdfstring{$A_5(8,6)$}{A5(8,6)}.}\label{sec586}

 In this section we consider the case~$(n,d)_q=(8,6)_5$ from Example~$\ref{586}$. Corollary~\ref{1mod4} implies that~$A_5(8,6) \leq 70$. Using small computer experiments, we will obtain~$A_5(8,6) \leq 65$. 
 
  As in the proof of Theorem~$\ref{importantth}$, we will compare upper and lower bounds on~$|X|$. But since an~$(8,6)_5$-code~$C$ of size at most~$ 70$ does not necessarily contain a~$15$-block (as~$70=5 \cdot 14$), we need information about~$14$-blocks. To this end we show, using an analogous approach as in~$\cite{bogzin}$ (based on occurrences of symbols in columns of an equidistant code):

\begin{proposition}\label{14block}
Any~$(7,6)_5$-code~$C$ of size~$14$ can be extended to a~$(7,6)_5$-code of size~$15$.
\end{proposition}
\proof 
For~$M=14$, the leftmost term in~$(\ref{detruc})$ equals the rightmost term. So~$C$ is equidistant and for each~$j \in \{1,\ldots,7\}$ there exists a unique~$\beta_j \in [q]$ with~$c_{\beta_j,j} =2$ and~$c_{\alpha,j}=3$ for all~$\alpha \in [q] \setminus \{\beta_j \}$. We can define a~$15$-th codeword~$u$ by putting~$u_j:= \beta_j$ for all~$j=1,\ldots,7$. We claim that~$C \cup \{ u \}$ is a~$(7,6)_5$-code of size~$15$.

To establish the claim we must prove that~$d_H(u,w) \geq 6$ for all~$w \in C$. Suppose that there is a word~$w \in C$ with~$d_H(u,w) <6$. We can renumber the symbols in each column of~$C$ such that~$w=\mathbf{1}$. Since~$C$ is equidistant, each word in~$C \setminus \{ w \}$ contains precisely one~$1$. On the other hand, there are two column indices~$j_1$ and~$j_2$ with~$u_{j_1} = 1$ and~$u_{j_2} = 1$.  Then~$C\setminus\{w\}$ contains at most~$1+1+5\cdot 2 = 12$ occurrences of the symbol~$1$ (since in columns~$j_1$ and~$j_2$ there is precisely one~$1$ in~$C \setminus \{w\}$). But in that case, since~$|C\setminus \{ w \}|=13>12$, there is a row in~$C$ that contains zero occurrences of the symbol~$1$, contradicting the fact that~$C$ is equidistant.
\endproof 

\noindent Note that a code of size more than~$65$ must have at least one~15- or~$14$-block, and therefore it must have a subcode of size~$65$ containing at least one~$15$- or~$14$-block. We shall now prove that this is impossible because 
\begin{align} \label{only13}
\text{each~$(8,6)_5$-code of size~$65$ only admits~$13$-blocks.} 
\end{align} 
It follows that~$A_5(8,6) \leq 65$. In order to prove~$(\ref{only13})$, let~$C$ be a~$(8,6)_5$-code of size~$65$. We first compute a lower bound on the number of irregular pairs in~$C$. Define, for~$x,y \in \mathbb{Z}_{\geq 0}$,
\begin{align}\label{lowX}
    f(x,y) &:= (3 x+ y)(65-15  x-14 y) + 3 \cdot 15 \binom{x}{2} +  14  \binom{y}{2} + 3\cdot 14 xy
    \\& \phantom{={,}}  - 2 \cdot 21  x- 8 y + \mathbf{1}_{\{y>0 \text{ and } x=0\}}  (65-14-39). \notag 
\end{align}

\begin{proposition}[Lower bound on~$|X|$]\label{LBX} Let~$C$ be an~$(n,d)_q=(8,6)_5$-code of size~$65$ and let~$j \in [n]$. Let~$x$ and~$y$ be the number of symbols that appear~$15$ and~$14$ times (respectively) in column~$j$. Then the number~$|X|$ of irregular pairs in~$C$ is at least~$f(x,y)$.
\end{proposition} 
\proof 
First consider a~$(7,6)_5$-code~$D$ of size~$15$ or size~$14$ and define 
\begin{align} 
S := \{ u \in [5]^7 \,\, | \,\,d_H(w,u) \geq 5 \,\,\,\, \forall\, w \in D \}. 
\end{align} 
For any~$u \in S$, define
\begin{align}
\alpha(u) := |\{w \in D\,\,: \, \, d_H(u,w)=6 \}|.
\end{align}
Then
\begin{align} \label{calign}
\text{if $|D|=15$, then} \phantom{aaai}   &      &   \text{if $|D|=14$, then}\phantom{aaai}&  \\
| \{u \in S\,\, |\,\, \alpha(u)=0\}| &=0,  &  | \{u \in S\,\, |\,\, \alpha(u)=0\}| &\leq 8, \notag \\
| \{u \in S\,\, |\,\, \alpha(u)=1\}| &\leq 21,  &  |\{ u \in S\,\, |\,\, \alpha(u)\leq1\}| &\leq 39. \notag \\
| \{u \in S\,\, |\,\, \alpha(u)=2\}|&=0.    &  &\notag 
\end{align}
This can be checked efficiently with a computer\footnote{All computer tests in this paper are small and can be executed within a minute on modern personal computers.} by checking all possible~$(7,6)_5$-codes  of size~$15$ and~$14$ up to equivalence. Here we note that a~$(7,6)_5$-code~$D$ (which must be equidistant, see Example~\ref{586}) of size~$15$ corresponds to a solution to \emph{Kirkman's school girl problem}~$\cite{zinoviev}$.\footnote{Kirkman's school girl problem asks to arrange 15 girls 7 days in a row in groups of 3 such that no two girls appear in the same group twice. The 1-1-correspondence between~$(n,d)_q=(7,6)_5$-codes~$D$ of size~$15$ and solutions to Kirkman's school girl problem is given by the rule: $\text{girls~$i_1$ and~$i_2$ walk in the same triple on day~$j$ } \Longleftrightarrow D_{i_1,j} = D_{i_2,j}$. }
 So to establish~$(\ref{calign})$, it suffices to check\footnote{By `check' we mean that given a~$(7,6)_5$-code~$D$ of size~$14$ or~$15$, we first compute~$S$, then~$\alpha(u)$ for all~$u \in S$, and subsequently verify~$(\ref{calign})$.} all~$(7,6)_5$-codes of size~$15$, that is, Kirkman systems (there are 7 nonisomorphic Kirkman systems~$\cite{7sol}$), and all~$(7,6)_5$-codes of size~$14$, of which there are at most~$7 \cdot 15$ by Proposition~$\ref{14block}$. 
 
Let~$G=(C,X)$ be the graph with vertex set~$V(G):=C$ and edge set~$E(G):=X$. Consider a~$15$-block~$B$ determined by column~$j$.  By~$(\ref{calign})$, each~$u \in C \setminus B$ has~$\geq 1$ neighbour in~$B$. We observed this also in Example~$\ref{586}$: for any~$u \in C \setminus B$ there exists at least one~$v \in B$ such that~$d_H(u,v) \notin \{6,8\}$, so~$d_H(u,v)=7$ and~$\{u,v\} \in X$. In~$(\ref{calign})$ this is represented as: if~$|D|=15$ then~$|\{u \in S\,\, | \,\, \alpha(u)=0\}|=0$, i.e., for any word~$u'$ of length~$7$ that has distance~$\geq 5$ to all words in a~$(7,6)_5$-code~$D$ of size~$15$, there is at least one~$v'\in D$ such that~$d_H(u',v')=6$. 

Furthermore,~$(\ref{calign})$ gives that all but~$\leq 21$ elements~$u\in C\setminus B$ have~$\geq 3$ neighbours in~$B$. So by adding~$\leq 2 \cdot 21$ new edges, we obtain that each~$u \in C \setminus B$ has~$\geq 3$ neighbours in~$B$. 

Similarly, for any~$14$-block~$B$ determined by column~$j$, by adding~$\leq 8$ new edges we achieve that each~$u \in C \setminus B$ has~$\geq 1$ neighbour in~$B$. Hence, by adding~$\leq ( 2 \cdot 21 \cdot x+  8\cdot y )$ edges to~$G$, we obtain a graph~$G'$ with
\begin{align}
|E(G')| \geq (3 x+ y)(65-15 x-14 y) + 3 \cdot 15  \binom{x}{2} +  14 \binom{y}{2}+ 3\cdot 14  xy.
\end{align}
 This results in the required bound, except for the term with the indicator function. That term can be added because $|\{ u \in S\,\, |\,\, \alpha(u)\leq1\}| \leq 39$ if~$|D|=14$, by~$(\ref{calign})$.
\endproof

\noindent It is also possible to give an upper bound on~$|X|$. If~$D$ is a~$(7,6)_5$-code of size~$k$, an upper bound~$h(k)=L-R$ on the number of pairs~$\{u,v\} \subseteq D$ with~$u \neq v$ and~$d_H(u,v)\neq 6$ (hence~$d_H(u,v)=7$) is given by~$(\ref{boundnd})$.  The resulting values~$h(k)$ are given in Table~\ref{UB}.

\begin{table}[ht]
\centering
\begin{tabular}{l|lllllllllll}
$k$ & 15 & 14 & 13 & 12 &11  &10  &9 & 8 & 7 & 6 & 5  \\\hline
$h(k)$      &0  & 0 & 1 & 3 & 6 & 10  & 8 & 7 & 7 & 8 & 10 
\end{tabular}
\caption{Upper bound~$h(k)$ on the number of pairs~$\{u,v\} \subseteq D$ with~$d_H(u,v)=7$ for a~$(7,6)_5$-code~$D$ with~$|D|=k$.}
\label{UB}
\end{table}

\begin{theorem}[$A_{5}(8,6) \leq 65$]\label{65th} 
Suppose that~$C$ is an~$(n,d)_q=(8,6)_5$-code with~$|C|=65$. Then each symbol appears exactly~$13$ times in each column of~$C$. Hence,~$A_{5}(8,6) \leq 65$.
\end{theorem}
\proof 
Let~$a^{(j)}_{k}$ be the number of symbols that appear exactly~$k$ times in column~$j$ of~$C$. Then the number of irregular pairs that have the same entry in column~$j$ is at most~$\sum_{k =5}^{15} a^{(j)}_{k}h(k)$. It follows that 
\begin{align}\label{contra}
    |X| \leq U := \sum_{j=1}^8 \sum_{k =5}^{15} a^{(j)}_{k} h(k).
\end{align}
One may check that if~$\mathbf{a}, \mathbf{b} \in \mathbb{Z}_{\geq 0}^{15}$ are $15$-tuples of nonnegative integers, with $\sum_k a_k k =65$, $\sum_k b_k k =65$, $\sum_k a_k=5$, $\sum_k b_k=5$, and $f(a_{15},a_{14}) \leq f(b_{15},b_{14}) \neq 0$, then
\begin{align} \label{tocheck}
\sum_{k=5}^{15} (7a_k +b_k) h(k) < f(b_{15},b_{14}).
\end{align}
(There are~$30$~$\mathbf{a} \in \mathbb{Z}_{\geq 0}^{15}$ with~$\sum_k a_k k =65$  and~$\sum_{k}a_k = 5$. So there are~$900$ possible pairs~$\mathbf{a},\mathbf{b}$. A computer now quickly verifies~$(\ref{tocheck})$.)    

By permuting the columns of~$C$ we may assume that $\max_j f( a^{(j)}_{15},a^{(j)}_{14})=f(a^{(1)}_{15},a^{(1)}_{14})$. Hence if~$f(a^{(1)}_{15},a^{(1)}_{14})>0$, then
\begin{align}
U &= \sum_{j=1}^8 \sum_{k =5}^{15} a^{(j)}_{k} h(k) = \frac{1}{7}\sum_{j=2}^8\left( \sum_{k =5}^{15} \left(7a^{(j)}_{k}  + a^{(1)}_{k} \right) h(k) \right) 
\\&< f(a^{(1)}_{15},a^{(1)}_{14}) \leq |X| \notag 
\end{align}
(where we used Proposition~$\ref{LBX}$ in the last inequality), contradicting~$(\ref{contra})$. So~$f( a^{(j)}_{15},a^{(j)}_{14})=0$ for all~$j$, which implies (for~$\mathbf{a}^{(j)} \in \mathbb{Z}_{\geq 0}^{15}$ with~$\sum_k a^{(j)}_k k =65$, $\sum_k a^{(j)}_k=5$) that~$a^{(j)}_{15}=a^{(j)}_{14}=0$ for all~$j$, hence each symbol appears exactly~$13$ times in each column of~$C$.
\endproof 

\begin{corollary}
$A_5(9,6) \leq 325$,~$A_5(10,6) \leq 1625$ and~$A_5(11,6) \leq 8125$. 
\end{corollary}
\proof 
By Theorem~$\ref{65th}$ and~$(\ref{elementarybounds2})$.
\endproof 

\section{Improved bound on~\texorpdfstring{$A_3(16,11)$}{A3(16,11)}.}\label{4118}

We show that~$A_3(16,11) \leq 29$ using a surprisingly simple argument.  

\begin{proposition}
$A_3(16,11) \leq 29$. 
\end{proposition}
\proof 
Suppose that~$C$ is an~$(n,d)_q=(16,11)_3$-code of size~$30$. We can assume that~$\mathbf{1} \in C$. It is known that~$A_3(15,11)=10$, so the symbol~$1$ is contained at most~$10$ times in every column of~$C$. Since~$|C|=30$, the symbol~$1$ appears exactly~10 times in every column of~$C$, so the number of 1's in~$C$ is divisible by~$5$. On the other hand it is easy to check  that a~$(15,11)_3$-code of size~$10$ is equidistant (using~$(\ref{boundnd})$, as~$L=R$). This implies that all distances in a~$(16,11)_3$-code of size~$30$ belong to~$\{11, 16\}$. So the number of~$1$'s in any code word~$\neq \mathbf{1}$ is~$0$ or~$5$. As~$\mathbf{1}$ contains~$16$ 1's, it follows that the total number of~1's is not divisible by~$5$, a contradiction.
\endproof

\section{Codes from symmetric nets}\label{symsec}

In this section we will show that there is a~$1$-$1$-correspondence between \emph{symmetric $(\mu ,q)$-nets} and $(n,d)_q=(\mu q,\mu q-\mu)_q$-codes of size~$\mu q^2$. From this, we derive in Section~$\ref{496}$ the new upper bound~$A_4(9,6) \leq 120$, implying~$A_4(10,6) \leq 480$.

\begin{defn}[Symmetric net] Let~$\mu,q \in \mathbb{N}$. A \emph{symmetric $(\mu, q)$-net} (also called \emph{symmetric transversal design}~$\cite{beth}$) is a set~$X$ of~$\mu q^2$ elements, called \emph{points}, together with a collection~$\mathcal{B}$ of subsets of~$X$ of size~$\mu q$, called \emph{blocks}, such that:
\begin{enumerate}[(s1)]
  \setlength{\itemsep}{1pt}
  \setlength{\parskip}{0pt}
  \setlength{\parsep}{0pt}
\item $\mathcal{B}$ can be partitioned into~$\mu q$ partitions (\emph{block parallel classes}) of~$X$.
\item Any two blocks that belong to different parallel classes intersect in exactly~$\mu$ points.
\item $X$ can be partitioned into~$\mu q$ sets of~$q$ points (\emph{point parallel classes}), such that any two points from different classes occur together in exactly~$\mu$ blocks, while any two points from the same class do not occur together in any block.\footnote{That is, a symmetric~$(\mu,q)$-net is a~$1-(\mu q^2, \mu q, \mu q)$ design~$D$, which is resolvable (s1), affine (s2), and the dual design~$D^*$ of~$D$ is affine resolvable (s3).}
\end{enumerate} 
\begin{remark}  From the 1-1-correspondence between symmetric~$(\mu,q)$-nets and $(n,d)_q=(\mu q,\mu q - \mu)_q$-codes~$C$ of size~$\mu q^2$ in Theorem~\ref{1rel} below it follows that~(s2) and~(s3) can be replaced by the single condition: 
\begin{itemize} 
\item[(s')] Each pair of points is contained in at most~$\mu$ blocks,
\end{itemize} 
since the only condition posed on such a code is that~$g(u,v) \leq \mu$ for all distinct~$u,v \in C$. 
\end{remark}
\end{defn}

\begin{examp}
Let~$X=\{1,2,3,4\}$ and~$\mathcal{B}=\{\{1,3\},\,\{2,4\},\,\{1,4\},\,\{2,3\} \}$. Then~$(X,\mathcal{B})$ is a symmetric~$(1,2)$-net. The block parallel classes are~$\{\{1,3\},\,\{2,4\}\}$ and~$\{\{1,4\},\,\{2,3\}\}$. The point parallel classes are~$\{1,2\}$ and~$\{3,4\}$. 
\end{examp}

\noindent By labeling the points as~$x_1,\ldots,x_{\mu q^2}$ and the blocks as~$B_1,\ldots,B_{\mu q^2}$,  the~\emph{$\mu q^2 \times \mu q^2$-incidence matrix}~$N$ of a symmetric $(\mu,q)$-net is defined by
\begin{align}
N_{i,j} := \begin{cases}1 &\mbox{if } x_i \in B_j, \\
0 &\mbox{else}.
\end{cases}
\end{align}
\noindent An \emph{isomorphism} of symmetric nets is a bijection from one symmetric net to another symmetric net that maps the blocks of the first net into the blocks of the second net. That is, two symmetric nets are isomorphic if and only if their incidence matrices are the same up to row and column permutations. Symmetric nets are, in some sense, a generalization of \emph{generalized Hadamard matrices}. 

\begin{defn}[Generalized Hadamard matrix]
Let~$M$ be an~$n \times n$-matrix with entries from a finite group~$G$. Then~$M$ is called a \emph{generalized Hadamard matrix} GH$(n,G)$ (or~GH$(n,|G|)$) if for any two different rows~$i$ and~$k$, the~$n$-tuple $(M_{ij}M_{jk}^{-1})_{j=1}^n$ contains each element of~$G$ exactly~$n/|G|$ times. 
\end{defn}

\begin{figure}[ht]
\begin{align*} 
\left( \begin{array}{cccccccc}
e & e &e &e &e &e &e & e\\
e & e &  a &a &b &b &c &c  \\
e & b & e  &b & c& a&c &a  \\
e &  c&  c &e &a & b& b& a \\
e & a & b  &c &e &a & b& c \\
e &  c&  b & a&c &e &a & b \\
e &b  & a  & c& a& c&e & b \\
e &  a&  c &b &b &c & a&e  \\ \end{array} \right)
\end{align*}
    \caption{An incidence matrix of the unique (up to isomorphism) symmetric~$(2,4)$-net is obtained by writing the elements~$e,a,b,c$ as $4 \times 4$-permutation matrices in the generalized Hadamard matrix GH$(8,V_4)$ (with~$V_4$ the Klein 4-group). See Al-Kenani~$\cite{42net}$.}
\end{figure} 

\noindent Each generalized Hadamard matrix~GH$(n,G)$ gives rise to a symmetric~$(n/|G|,|G|)$-net: by replacing~$G$ by a set of~$|G| \times |G|$-permutation matrices isomorphic to~$G$ (as a group), one obtains the incidence matrix of a symmetric net. Not every symmetric~$(n/q,q)$-net gives rise to a generalized Hadamard matrix~GH$(n,q)$, see~$\cite{matrixnet}$. But if the group of automorphisms (\emph{bitranslations}) of a symmetric $(n/q,q)$-net has order~$q$, then one can construct a generalized Hadamard matrix~GH$(n,q)$ from it. See~$\cite{beth}$ for details. 

\begin{assumption} 
In this section we consider triples~$(n,d)_q$ of natural numbers for which
\begin{align}
qd = (q-1)n,
\end{align} 
hence~$n-d=n/q=:\mu$ and~$\mu \in \N$. So~$(n,d)_q=(\mu q, \mu q - \mu)_q$. 
\end{assumption}

\noindent The fact that a generalized Hadamard matrix~$\text{GH}(n,q)$ gives rise to an~$(n,d)_q$-code of size~$qn$, was proved in~$\cite{plotkin}$ and for some parameters it can also be deduced from an earlier paper~$\cite{zinoviev2}$. Using a result by Bassalygo, Dodunekov, Zinoviev and Helleseth~$\cite{granking}$ about the structure of~$(n,d)_q$-codes of size~$qn$,\footnote{Note that~$A_q(n,d)\leq qn$, since by Plotkin's bound~$(\ref{elementarybounds})$, $A_q(n-1,d) \leq n$, hence~$A_q(n,d) \leq qn =\mu q^2$ by~$(\ref{elementarybounds2})$.} we prove that such codes are in 1-1-relation with symmetric~$(n/q,q)$-nets. 

\begin{theorem} \label{1rel} Let~$\mu, q \in \N$. There is a~1-1-relation between  symmetric~$(\mu,q)$-nets (up to isomorphism) and $(n,d)_q=(\mu q,\mu q -\mu)_q $-codes~$C$ of size~$\mu q^2$ (up to equivalence). 
\end{theorem}
\proof 
Given an $(n,d)_q= (\mu q,\mu q -\mu)_q$-code~$C$ of size~$\mu q^2$, we construct a $(0,1)$-matrix~$M$ of order $\mu q^2 \times \mu q^2$ with the following properties:
\begin{enumerate}[(I)]
\item \label{I} $M$ is a~$\mu q^2 \times \mu q^2$ matrix that consists of~$q \times q$ blocks~$\sigma_{i,j}$ (so~$M$ is a~$\mu q \times \mu q$ matrix of blocks~$\sigma_{i,j}$), where each~$\sigma_{i,j}$ is a permutation matrix.
\item \label{II} $MM^T= M^TM= A$, where~$A$ is a~$\mu q^2 \times \mu q^2$ matrix that consists of~$q \times q$ blocks~$A_{i,j}$ (so~$A$ is an~$\mu q \times \mu q$ matrix of blocks~$A_{i,j}$), with
\begin{align} \label{M}
A_{i,j} =\begin{cases} \mu q \cdot I_q &\mbox{if } i =j, \\
\mu \cdot J_q & \mbox{if } i \neq j. \end{cases}
\end{align}
Here~$J_q$ denotes the~$q \times q$ all-ones matrix.
\end{enumerate}
 By Proposition 4 of~$\cite{granking}$, since~$d=n(q-1)/q$ and~$|C|=qn$, $C$ can be partitioned as
\begin{align} \label{partition} 
C = V_1 \cup V_2 \cup \ldots \cup V_{n},
\end{align}
where the union is disjoint,~$|V_i|=q$ for all~$i=1,\ldots,n$, and where~$d_H(u,v)=n$ if~$u,v \in C$ are together in one of the~$V_i$, and~$d_H(u,v)=d$ if~$u\in V_i$ and~$v \in V_j$ with~$i\neq j$.

Now we write each word~$w \in [q]^{n}$ as a $(0,1)$-row vector of size~$qn = \mu q^2$ by putting a~$1$ on positions~$(i,w_i) \in [n] \times [q]$ (for~$i=1,\ldots,n$) and 0's elsewhere. The $q$~words in any of the~$V_i$ then form a~$q \times qn$ matrix consisting of~$n$ permutation matrices~$\sigma_{i,j}$ of size~$q \times q$. 

By placing the matrices obtained in this way from all~$n$ tuples~$V_1,\ldots,V_{n}$  underneath each other, we obtain a~$qn \times qn$ matrix~$M$ consisting of~$n^2$ permutation matrices of order~$q\times q$, so~(\ref{I}) is satisfied. Property~(\ref{II}) also holds, since for any~$u,v \in C$ written as row vectors of size~$qn$, with the~$V_i$ as in~$(\ref{partition})$, it holds that
\begin{align} 
\sum_{k \in  [n] \times [q]}u_kv_k= g(u,v)=\begin{cases}n = \mu q &\mbox{if } u = v, \\
0  &\mbox{if } u \neq v \text{ and } u,v \in V_i,\\
n-d=\mu  &\mbox{if } u \neq v \text{ and } u \in V_i, v \in V_j \text{ with } i \neq j.\\
\end{cases}
\end{align}
So~$MM^T=A$. Moreover, if~$j_1:=(j_1',a_1) \in [n] \times [q]$ and~$j_2:=(j_2',a_2) \in [n] \times [q]$, then
\begin{align} 
\sum_{k \in  [qn]}M_{k,j_1} M_{k,j_2}=\begin{cases}n = \mu q  &\mbox{if } j_1'=j_2' \text{ and } a_1=a_2, \\
0  &\mbox{if } j_1' = j_2' \text{ and }  a_1\neq a_2,\\
n/q = \mu  &\mbox{if } j_1' \neq j_2',\\
\end{cases}
\end{align}
where the last statement follows by considering the words in~$C$ that have~$a_1$ at the~$j_1'$-th position. (The remaining columns form an~$n$-block for the~$j_1'$-th column. In this~$n$-block, each symbol occurs exactly~$n/q$ times at each position, since~the leftmost term equals the rightmost term in~$(\ref{detruc})$ for~$(n-1,d)_q$-codes of size~$n$.) We see that also~$M^TM=A$. Hence,~$M$ is the incidence matrix of a symmetric~$(\mu,q)$-net (see~$\cite{beth}$, Proposition~I.7.6 for the net and its dual). 

Note that one can do the reverse construction as well: given a symmetric~$(\mu,q)$-net, the incidence matrix of~$M$ can be written (after possible row and column permutations) as a matrix of permutation matrices such that~$MM^T=M^TM=A$, with~$A$ as in~$(\ref{M})$. From~$M$ we obtain a code~$C$ of size~$\mu q^2$ of the required minimum distance by mapping the rows $(i,w_i) \in [\mu q] \times [q]$ to~$w \in [q]^{\mu q}$. Observe that equivalent codes yield isomorphic incidence matrices~$M$ and vice versa.
\endproof 

\begin{figure}[H]
{\footnotesize$$
\begin{blockarray}{c|ccc|}
    & & &  \\  \cline{1-4}
 \begin{block}{c|ccc|}
   w_1 &   1 & 1 & 1 \\ 
   w_2 &   2 & 2 & 2 \\
   w_3 &   3 & 3 & 3 \\ 
   w_4 &   1 & 3 & 2 \\ 
   w_5 &   2 & 1 & 3 \\ 
   w_6 &   3 & 2 & 1 \\ 
   w_7 &   1 & 2 & 3 \\ 
   w_8 &   2 & 3 & 1 \\  
   w_9 &   3 & 1 & 2 \\    \end{block}
\end{blockarray} \,\, \quad \longleftrightarrow
\,\,\quad 
\begin{blockarray}{c|ccc|ccc|ccc|}
 &  1 & 2 & 3 & 1 & 2 & 3 & 1 & 2 & 3  \\  \cline{1-10}
 \begin{block}{c|ccc|ccc|ccc|}
   w_1 &   1 & 0 & 0 & 1 & 0 &0  & 1 & 0 & 0\\
   w_2 &   0 & 1 & 0 & 0 & 1 & 0 & 0 & 1 & 0\\
   w_3 &   0 & 0 & 1 & 0 & 0 & 1 & 0 & 0 & 1\\\cline{1-10}
   w_4 &   1 &0  & 0 & 0 & 0 & 1 & 0 & 1 &0 \\
   w_5 &   0 & 1 & 0 & 1 & 0 & 0 & 0 & 0 &1 \\
   w_6 &   0 &  0&  1& 0 & 1 & 0 & 1 & 0 &0 \\\cline{1-10}
   w_7 &   1 & 0 &  0& 0 & 1 & 0 & 0 & 0 &1 \\
   w_8 &   0 & 1 &  0& 0 & 0 & 1 & 1 & 0 & 0\\
   w_9 &   0 & 0 & 1 & 1 & 0 & 0 & 0 & 1 &0 \\  \end{block}
\end{blockarray} 
$$ }
    \caption{An~$(n,d)_q=(3,2)_3$-code~$C=\{w_1,\ldots,w_9\}$ of size~$9$ (left table) gives rise to an incidence matrix of a symmetric~$(1,3)$-net (right table) and vice versa.}
\end{figure}

\section{New upper bound on~\texorpdfstring{$A_4(9,6)$}{A4(9,6)}.}\label{496}

\noindent In this section we use the 1-1-correspondence between symmetric~$(\mu,q)$-nets and~$(n,d)_q=(\mu q, \mu q - \mu)_q$-codes of size~$\mu q^2$ in combination with a known result about symmetric~$(2,4)$-nets~$\cite{42net}$ to derive that~$A_4(9,6) \leq 120$.

As~$A_4(8,6)=32$, any~$(9,6)_4$-code of size more than~$120$ must contain at least one~$31$- or~$32$-block, and therefore it contains a subcode of size~$120$ containing at least one~$31$- or~$32$-block. We will show (using a small computer check) that this is impossible because a~$(9,6)_4$-code of size~$120$ does not contain any~$31$- or~$32$-blocks. Therefore~$A_4(9,6) \leq 120$.  In order to do prove this, we need information about~$(8,6)_4$-codes of size~$31$.

\begin{proposition}\label{31code}
Let~$q,n,d \in \N$ satisfy~$qd=(q-1)n$. Any~$(n,d)_q$-code~$C$ of size~$qn-1$ can be extended to an~$(n,d)_q$-code of size~$qn$.
\end{proposition}
\proof 
Let~$C$ be an~$(n,d)_q$-code of size~$qn-1$. By Plotkin's bound, $A_q(n-1,d) \leq n$, so each symbol occurs at most~$n$ times in each column of~$C$, hence there exists for each~$j \in [n]$ a unique~$\beta_j \in [q]$ with~$c_{\beta_j,j} =n-1$ and~$c_{\alpha,j}=n$ for all~$\alpha \in [q] \setminus \{\beta_j \}$. We can define a~$qn$-th codeword~$u$ by putting~$u_j:= \beta_j$ for all~$j=1,\ldots,n$. We claim that~$C \cup \{ u \}$ is an~$(n,d)_q$-code of size~$qn$.

To establish the claim we must prove that~$d_H(u,w) \geq d$ for all~$w \in C$. Let~$w \in C$ with~$d_H(u,w) < n$. We can renumber the symbols in each column of~$C$ such that~$w=\mathbf{1}$. Then~$w$ is contained in an~$(n-1)$-block~$B$ for some column in~$C$ (otherwise~$d_H(u,w)=n$). The number of~1's in~$B$ is~$n+(n-2)n/q$ (since any $(q,n-1,d)$-code of size~$n-1$ is equidistant, as~$L-R=0$ in~$(\ref{boundnd})$ for~$(n-1,d)_q$-codes of size~$n-1$) and the number of~1's in~$C \setminus B$ is~$(q-1)(n-1)n/q$ (since in any $(n-1,d)_q$-code of size~$n$, each symbol appears exactly~$n/q$ times in each column, as the leftmost term equals the rightmost term in~$(\ref{detruc})$ for~$(n-1,d)_q$-codes of size~$n$). Adding these two numbers we see that the number of~1's in~$C$ is~$n^2-n/q$. Since~$C \cup \{u\}$ contains each symbol~$n^2$ times by construction,~$u$ contains symbol~$1$ exactly~$n/q$ times, hence~$d_H(u,w)=n-n/q=d$, which gives the desired result.
\endproof

\begin{proposition}\label{120}
$A_4(9,6) \leq 120$. 
\end{proposition} 
\proof 
The~$(n,d)_q=(8,6)_4$-code of size~$32$ is unique up to equivalence, since the symmetric~$(2,4)$-net is unique up to equivalence (see Al-Kenani~$\cite{42net}$). By checking all~$(8,6)_4$-codes~$D$ of size~$31$ (of which there are at most~$32$ up to equivalence since each~$(8,6)_4$-code of size~$31$ arises by removing one word from a~$(8,6)_4$-code of size~$32$ by Proposition~\ref{31code}) we find that
\begin{align} 
|\{ u \in [4]^{8} \,\, | \,\,d_H(w,u) \geq 5 \,\,\,\, \forall\, w \in D \}| \leq 25. 
\end{align} 
This implies that an $(n,d)_q=(9,6)_4$-code~$C$ of size~$120$ cannot contain a~$31$- or~$32$-block. Therefore $A_4(9,6) \leq 120$.
\endproof 

\begin{corollary}
$A_4(10,6) \leq  480$.
\end{corollary}
\proof 
By Proposition~$\ref{120}$ and~$(\ref{elementarybounds2})$.
\endproof 

\section*{Acknowledgements}
I am most grateful to Lex Schrijver for his supervision, for his help regarding both the content and the presentation of the paper and for all the conversations we had about the code bounds (in person and by e-mail). Thank you, Lex!

Also I am grateful to Bart Litjens, Guus Regts and Jacob Turner for their comments. Furtermore I want to thank the editor and the anonymous referees for their very helpful comments concerning the presentation of the material.

\selectlanguage{english}

\end{document}